\documentclass[a4paper,10pt]{amsart}
\usepackage{amscd}

\def\oM{\overline{M}}

\def\bP{\Bbb P}

\def\Hilb{\operatorname{Hilb}}

\def\cO{\Cal O}
\def\cN{\Cal N}

\def\Bl{\operatorname{Bl}}

\def\bZ{\Bbb Z}

\def\bQ{\Bbb Q}
\def\bG{\Bbb G}
\def\bR{\Bbb R}
\def\bA{\Bbb A}

\def\oX{\overline{X}}

\def\cB{\Cal B}
\def\cA{\Cal A}

\def\cF{\Cal F}

\def\cI{\Cal I}
\def\cA{\Cal A}

\def\cB{\Cal B}

\def\cB{\Cal B}

\def\cH{\Cal H}

\def\Spec{\operatorname{Spec}}

\def\cK{\Cal K}

\def\cO{\Cal O}

\def\oK{\overline{K}}

\def\Spec{\operatorname{Spec}}
\def\pr{\operatorname{pr}}
\def\Hom{\operatorname{Hom}}
\def\Cal{\mathcal}

\newtheoremstyle{mystyle}{}{}{\itshape}{}{\scshape}{.}{ }{}
\theoremstyle{mystyle}

\swapnumbers
\newtheorem{Theorem}{Theorem}[section]

\newtheorem{Proposition}[Theorem]{Proposition}
\newtheorem{Lemma}[Theorem]{Lemma}

\newtheoremstyle{myreview}{}{}{}{}{\scshape}{.}{ }{}
\theoremstyle{myreview}

\newtheorem{Definition}[Theorem]{Definition}
\newtheorem{Example}[Theorem]{Example}

\newtheorem{Remark}[Theorem]{Remark}

\newtheorem{Review}[Theorem]{}

\newcounter{et}[Theorem]
\def\cooltag{\tag{\arabic{section}.\arabic{Theorem}.\arabic{et}}\addtocounter{et}{1}}

\begin{document}

\title{Compactifications of Subvarieties of Tori}

\author{Jenia Tevelev}

\begin{abstract}
We study compactifications of subvarieties of algebraic tori
defined by imposing a sufficiently fine polyhedral structure
on their non-archimedean amoebas.
These compactifications have many nice properties,
for example any $k$ boundary divisors
intersect in codimension $k$.
We consider some examples including $M_{0,n}\subset\oM_{0,n}$
(and more generally log canonical models of complements of hyperplane arrangements)
and compact quotients of Grassmannians by a maximal torus.
\end{abstract}

\maketitle

\section{Introduction and Statement of Results}

Let $X$ be a connected closed subvariety of an algebraic torus $T$
over an algebraically closed field $k$. It is natural to consider
compactifications $\oX$ defined as closures of $X$ in various toric
varieties $\bP$ of $T$. In this paper we address the following
question: how nice could $\oX$ possibly be? Consider the
multiplication map
$$\Psi:\,T\times\oX\to\bP, \quad (t,x)\mapsto tx.$$

\begin{Definition}
We call $\oX$ a {\em tropical compactification} if $\Psi$ is
faithfully flat and $\oX$ is proper.
A term ``tropical'' is
explained in Section~\ref{main}.
\end{Definition}

\begin{Theorem}\label{MainTh}
Any subvariety $X$ of a torus has a tropical compactification~$\oX$
such that the corresponding toric variety $\bP$ is smooth. In this
case the boundary $\oX\setminus X$ is divisorial and has
``combinatorial normal crossings'': for any collection
$B_1,\ldots,B_r\in\oX\setminus X$ of irreducible divisors, $\cap
B_i$ has codimension $r$. Moreover, if $r=\dim X$ and $p\in
B_1\cap\ldots\cap B_r$ then $\oX$ is Cohen-Macaulay at~$p$.
\end{Theorem}

We also show that if $\oX\subset\bP$ is tropical then
$\oX'\subset\bP'$ is also tropical for any proper toric morphism
$\bP'\to\bP$. However, it is not known (at least to the author) if
there exists a minimal tropical compactification.

The boundary $\oX\setminus X$ does not necessarily have genuine
normal crossings.

\begin{Definition}
A subvariety $X$ of a torus $T$ is called {\em sch\"on} if it has a
tropical compactification with a smooth multiplication map.
\end{Definition}

\begin{Theorem}\label{KLHGKFJGFhgghhg}
If $X$ is sch\"on then any of its tropical compactifications
$\oX\subset\bP$ has a smooth multiplication map, is regularly
embedded, normal, has toroidal singularities, and its log canonical
line bundle is globally generated and is equal to the determinant of
the normal bundle. Moreover, if $\bP'\to\bP$ is a proper toric
morphism then the corresponding map $\oX'\to\oX$ is log crepant.
\end{Theorem}

Sch\"on hypersurfaces in tori are known as ``hypersurfaces
non-degenerate with respect to their Newton polytope'', introduced
and studied by Varchenko \cite{V1,V2}.

Next we discuss one natural class of sch\"on subvarieties of tori of
high codimension. Let $X$ be a complement of a hyperplane
arrangement,
\begin{equation}\label{KContour}
X=\bP^r\setminus\{H_1,\ldots,H_n\}.\tag{$\circ$}
\end{equation}
Assume that the arrangement is essential (not all hyperplanes pass
through a point) and connected (does not have any projective
automorphisms). $\bP^r$ admits an embedding in $\bP^{n-1}$ such that
the hyperplanes $H_i$ are intersections with coordinate hyperplanes
of~$\bP^{n-1}$. It follows that $X$ is a closed subvariety of a
torus $(k^*)^{n-1}$.

\begin{Theorem}\label{:LKJHGKGKJYF}\label{Hyper}\label{LKJGKFKJGFBVCH}
 $X$ is
sch\"on and Kapranov's visible contour $\oX_{vc}$ of $X$ \cite{CQ1}
is one of its tropical compactifications. The log canonical line
bundle of $\oX_{vc}$ is very ample and  $\oX_{vc}$ is the log
canonical model of $X$.
\end{Theorem}

This theorem was essentially proved in \cite[Theorem~2.19]{CQ2} but
we give another argument for the reader's convenience. We discuss
this result and compare the tropical compactification with other
known compactifications of complements of hyperplane arrangements in
Section~\ref{Hyperplane}.

In fact, our proof of existence of tropical compactifications is
constructive and blends the following ingredients:
\begin{itemize}
\item Generalized Kapranov's visible contours.
\item Lafforgue's transversality argument of \cite[Thm.~4.5]{Laf}.
\item Non-archimedean amoebas.
\end{itemize}

We will explain generalized visible contours below, everything else
is explained in Section~\ref{main}. Let $T_X$ be the normalizer of
$X$ in $T$. The action of $T_X$ on $T$ (and hence on $X$) is free
and the quotient $X/T_X$ is a closed subvariety of a torus $T/T_X$.
We will construct a tropical compactification of $X/T_X$. Since
$X\simeq (X/T_X)\times T_X$, it easily follows that $X$ also has a
tropical compactification.

\begin{Definition}\label{defvisiblecontour}
This construction depends on a choice of an auxiliary proper toric
variety $T\subset\bP$. Let $\oX$ be the closure of $X$ in $\bP$. For
any closed subscheme $Y$ of $\bP$, we denote by $[Y]\in\Hilb(\bP)$
the corresponding point of the Hilbert scheme. Consider the twisted
action of $T$ on $\Hilb(\bP)$:
$$T\times \Hilb(\bP)\to\Hilb(\bP),\quad t\cdot [Y]=t^{-1}[Y].$$

The normalization of the closure of the $T$-orbit of $[\oX]$ in
$\Hilb(\bP)$ is called a {\em Groebner toric variety} $\bP_{Gr}$ of
$(X,\bP)$. It compactifies the torus $T/T_X$. Let $\oX_{vc}$ be the
closure of $X/T_X$ in $\bP_{Gr}$. Let $\bP_{vc}\subset \bP_{Gr}$ be
the toric open subset consisting of all $T$-orbits that intersect
$\oX_{vc}$. We consider $\oX_{vc}$ as a closure of $X/T_X$ in
$\bP_{vc}$ (and not $\bP_{Gr}$) to make the multiplication map
surjective.
\end{Definition}

\begin{Theorem}\label{KJFKJFKD}\label{J<GKJFKJGFKHGF}
$\oX_{vc}$ is a tropical compactification of $X/T_X$. Moreover,
\begin{equation}\label{LKJGKJF}
\oX_{vc}=\{[Y]\in \bP_{Gr}\,|\,e\in Y\},\tag{$\star$}
\end{equation}
where $e\in T$ is the identity element.
\end{Theorem}

The formula \eqref{LKJGKJF} is a direct generalization of Kapranov's
visible contour \cite{CQ1}, which is a special case when $X$ is a
complement of a hyperplane arrangement as in~\ref{Hyper}, and
$\Hilb(\bP)$ is just the Grassmannian.

\begin{Remark}
Though $\oX_{trop}$ is defined using Hilbert schemes, in fact we use
only the closure of the $T$-orbit of $[\oX]$. This makes
$\oX_{trop}$ computable, see \ref{InitialIdeals}. A different method
of constructing tropical compactifications is discussed in
Remark~\ref{surjectiveamoeba}.
\end{Remark}

Any tropical compactification of $X/T_X$ can be viewed as a compact
quotient space of $X$. It turns out that it is related to (and
refines) other quotient constructions such as GIT quotients or more
precisely their inverse limits, Chow/Hilbert quotients. This
relationship is explained in Section~\ref{tropicalquotients}. We
concentrate in particular on tropical quotients of Grassmannians by
the maximal torus.

In what follows we assume that the reader is familiar with the
basics of the logarithmic Mori theory~\cite{KMM}.

\begin{Review}\textsc{Acknowledgements.}
I am grateful to Sean Keel for a lot of help, in particular he
suggested me Theorem~\ref{KLHGKFJGFhgghhg}. This paper is a part of
the project \cite{CQ2,HKT,KT} on Chow quotients of Grassmannians. I
am grateful to Paul Hacking, Misha Kapranov, Laurent Lafforgue,
Bernd Sturmfels, David Speyer, Martin Olsson, and David Saltman for
several illuminating discussions.
\end{Review}

\section{Tropical Compactifications}\label{main}
Let $X$ be a connected closed subvariety of an algebraic torus $T$
over an algebraically closed field $k$. Let $\oK$ be the field of
Puiseaux series over $k$ with the non-archimedean valuation
$\deg:\,\oK\to\bQ$ (the minimal exponent of the series). Let
$X(\oK)\subset T(\oK)$ be $\oK$-points. The valuation extends to the
surjective map
$$\deg:\,T(\oK)\to \Lambda_\bQ^\vee,$$
where $\Lambda^\vee$ is the $\Bbb Z$-lattice of 1-PS of $T$ and
$\Lambda_\bQ^\vee=\Lambda^\vee\otimes\bQ$.

\begin{Definition}[\cite{EKL}]
The {\em non-archimedean amoeba}  of a pair $(X,T)$ is
$$\cA=\deg X(\oK)\subset\Lambda_\bQ^\vee.$$
\end{Definition}

By \cite[Theorem 2.2.5]{EKL}, $\cA$ is the underlying set of a fan
(maybe not uniquely determined). More generally, if $X$ is defined
over $K$ then $\cA$ is the underlying set of a polyhedral complex
(maybe not uniquely determined) and is equal to the Bieri--Groves
set \cite{BG} of $X$.

\begin{Lemma}\label{KLUHKJG}
Let $n=\dim T$. Let $\cF\subset\Lambda_\bQ^\vee$ be a cone of
dimension $d$ such that its rays are spanned by vectors from a basis
of $\Lambda^\vee$. Let
$$\bP=\bA^d\times (k^*)^{n-d}$$ be the corresponding toric variety.
$\cF$ contains a point of $\cA$ in its interior if and only if the
closure $\oX$ of $X$ intersects the closed orbit of $\bP$.
\end{Lemma}

\begin{proof}
We can assume that $\Lambda^\vee$ has a basis $x_1,\ldots,x_n$
and that $\cF$ is spanned by $x_1,\ldots,x_d$.
We identify $T$ with $(k^*)^n$
and $\bP$ with $\bA^d\times(k^*)^{n-d}$  using this basis.
If $\cF$ contains a point of $\cA$ in its interior
then $X(\oK)$ contains a point
\begin{equation}\label{germ}
p(z)=(p_1(z),\ldots,p_n(z)),
\end{equation}
where $z$ is a uniformizer in~$K$, $\deg p_i(z)>0$ for $i\le d$, and $\deg p_i(z)=0$ for $i>d$.
It follows that $\oX(k)$ contains a point $p(0)=(0,\ldots,0,p_{d+1}(0),\ldots,p_n(0))$
that belongs to the closed orbit.
And vice versa, if $\oX$ intersects the closed orbit of $\bP$
then one can find a germ of a curve \eqref{germ} in $X$ such that
$\deg p_i(z)\ge 0$ for any $i$ and $(p_1(0),\ldots,p_n(0))$ belongs to the closed orbit.
In particular, $\deg p(z)$ belongs to the interior of $\cF$.
\end{proof}

\begin{Proposition}\label{comp}
Let $\bP$ be a toric variety of $T$ (not necessarily proper) with
the fan $\cF\subset\Lambda^\vee\otimes\bQ$. Let $\oX$ be the closure
of $X$ in $\bP$. Then $\oX$ is proper iff $\cF\supset\cA$.
\end{Proposition}

\begin{proof}
Suppose $\oX$ is proper but $\cA$ is not contained in $\cF$.
Consider any compactification $\bP\subset\bP'$ such that there
exists a simplicial cone $C$ of $\cF'$ with the following property.
Its interior $C^0$ does not intersect $\cF$ and contains a point of
$\cA$. Since $\oX$ is proper, we have $\oX=\oX'$, and therefore
$\oX$ does not intersect a $T$-orbit in $\bP'\setminus\bP$ that
corresponds to~$C$. This contradicts Lemma~\ref{KLUHKJG}.

Suppose $\oX$ is not proper but $\cA$ is contained in $\cF$.
Consider any compactification $\bP\subset\bP'$ such that all cones
in $\cF'\setminus\cF$ are simplicial. Since $\oX'$ is proper, it
intersects a $T$-orbit in $\bP'\setminus\bP$. But $\cA$ does not
intersect the interior of the corresponding cone of $\cF'$ which
again contradicts Lemma~\ref{KLUHKJG}.
\end{proof}

\begin{Definition}
We say that a pair $(X,\bP)$ is tropical if the closure $\oX$ of $X$
in $\bP$ is a tropical compactification of $X$, i.e. the
multiplication map $\Psi:\,T\times\oX\to\bP$, $(t,x)\mapsto tx$ is
faithfully flat and $\oX$ is proper.
\end{Definition}

\begin{Proposition}\label{:KLJGLJFKF}
Assume that $(X,\bP)$ is tropical. Then the fan of $\bP$ is
supported on $\cA$. Assume, moreover, that $\bP'\to\bP$ is any
proper toric morphism. Then $\oX'\subset\bP'$ is also tropical,
$\oX'$ is the inverse image of $\oX$, and $\Psi'$ is the pull back
of~$\Psi$.
\end{Proposition}

\begin{proof}
Assume $(X,\bP)$ is tropical and let $f:\,\bP'\to\bP$ be any proper
morphism. We claim that  $(X,\bP')$ is tropical. Since the pull-back
of a faithfully flat map is faithfully flat, it suffices to show
that
$$\Psi'=f^*\Psi,$$
i.e.~that $(\oX\times T)\times_{\bP}\bP'$ is integral. Since
$(X\times T)\times_{\bP}\bP'=X\times T$, this follows from
Lemma~\ref{jhgfdjdjfd} below. An obvious Cartesian diagram with
vertical arrows induced by~$f$
$$
\begin{CD}
\oX'  @>{(Id,e)}>> \oX'\times T @>{\Psi'=f^*\Psi}>> \bP'\\
@VVV  @VVV  @VVV \\
\oX  @>{(Id,e)}>> \oX\times T @>{\Psi}>> \bP\\
\end{CD}
$$
implies that $\oX'$ is an inverse image of $\oX$.

Finally, suppose that $(X,\bP)$ is tropical but $\cF$ is not
supported on $\cA$. Let $\cF'$ be a refinement of $\cF$ such that
$\cF'$ contains a cone $C$ such that its rays are spanned by vectors
from a basis of $\Lambda^\vee$ and such that $C^0\cap\cA=\emptyset$.
Then $(X,\bP')$ is tropical, and in particular $\Psi'$ is surjective
but by Lemma~\ref{KLUHKJG} $\oX'$ (and therefore the image of the
structure map) does not intersect a $T$-orbit in $\bP'$
corresponding to~$C$.
\end{proof}

\begin{Lemma}\label{jhgfdjdjfd}
Suppose $A\subset B$ is a flat extension of rings, $A$ is a domain, $f\in A$, and
the localization $B_f$ is a domain. Then $B$ is a domain.
\end{Lemma}

\begin{proof}
Since $A$ is a domain, $A\subset A_f$. Since $B$ is flat over~$A$,
$B\subset (A_f)\otimes_AB=B_f$. Therefore, $B$ is a domain.
\end{proof}

\begin{Remark}
We see that a tropical compactification is defined by a sufficiently
fine polyhedral structure on~$\cA$. It would be interesting to
refine Proposition~\ref{:KLJGLJFKF}, for example is it true that
$\oX$ is tropical if and only if $\cF$ is supported on~$\cA$? Is
there a canonical polyhedral structure on~$\cA$? Of course, we have
not yet proved that tropical compactifications even exist. This will
be proved later in this section.
\end{Remark}

\begin{proof}[Proof of Theorem~\ref{KLHGKFJGFhgghhg}]
Recall that a subvariety $X$ of a torus $T$ is called {\em sch\"on}
if it has a tropical compactification with a smooth multiplication
map. We claim that if $X$ is sch\"on then any of its tropical
compactifications $\oX\subset\bP$ has a smooth multiplication map,
is regularly embedded (i.e.~is locally defined by a regular
sequence), normal, has toroidal singularities, and its log canonical
line bundle is globally generated and is equal to the determinant of
the normal bundle of $\oX$. Moreover, if $\bP'\to\bP$ is a proper
toric morphism then the corresponding map $\oX'\to\oX$ is log
crepant.

Assume that $X$ is sch\"on and $(X,\bP_0)$ has a smooth structure
map. Let $(X,\bP)$ be another tropical object. We can find a
tropical object $(X,\bP')$ with morphisms $\bP'\to\bP$,
$\bP'\to\bP_0$. Then $\Psi'$ is the pullback of~$\Psi$ and the
pullback of $\Psi_0$ by Proposition~\ref{:KLJGLJFKF}. Therefore all
these maps have the same set of geometric fibers. But a flat map is
smooth if and only if all its geometric fibers are smooth. Thus
$\Psi$ is smooth too.

From the commutative diagram with smooth vertical arrows
$$
\begin{CD}
\oX\times T  @>{(\Psi,\pr_2)}>> \bP\times T\\
@V{\Psi}VV    @V{\pr_1}VV \\
\bP          & = &\bP\\
\end{CD}
$$
it follows that $(\Psi,\pr_2):\,\oX\times T\to \bP\times T$ is a
regular embedding by \cite[17.12.1]{EGA}. Since the map
$$\bP\times T\to\bP\times T,\quad (p,t)\mapsto(t^{-1}p,t)$$
is an automorphism, the obvious embedding $\oX\times
T\subset\bP\times T$ is also regular. This implies $\oX\subset\bP$
is a regular embedding.

Since $\Psi$ is smooth, $\oX\times T$ is normal and has toroidal
singularities. Therefore, $\oX$ has these properties as well.

Since $\oX\subset\bP$ is regularly embedded, it has a normal bundle
$\cN$. Since the log canonical bundle of a toric variety is trivial,
$\det\cN^\vee$ is the log canonical bundle of $\oX$ by adjunction.
Moreover,  for any toric morphism $\bP'\to\bP$, the $\cN'$ is the
pullback of $\cN$ because $\oX'$ is the inverse image of $\oX$ by
Proposition~\ref{:KLJGLJFKF}. Therefore, $\det(\cN')^\vee$ is the
pullback of $\det\cN^\vee$, i.e.~$\oX'\to\oX$ is log crepant.

Finally, from the canonical sequence
$$\Omega^1_{\oX}(\log)\to\Omega^1_{\bP}(\log)\to\cN^\vee\to0$$
of log tangent bundles we see that $\cN^\vee$ is globally generated
because the log tangent bundle of a toric variety is trivial.
Therefore, $\det\cN^\vee$ is also globally generated.
\end{proof}

Now we prove that tropical compactifications exist by proving that
the generalized visible contour compactification $\oX_{vc}$ is
tropical. We follow the notation of
Definition~\ref{defvisiblecontour} and Theorem~\ref{J<GKJFKJGFKHGF}.

\begin{proof}[Proof of Theorem~\ref{KJFKJFKD}]
We temporarily denote $\{[Y]\in \bP_{Gr}\,|\,e\in Y\}$ by $\cK$.
Then it is clear that $\oX_{vc}$ is contained in and is an
irreducible component of $\cK$: indeed, it suffices to show that
$x^{-1}\cdot \oX$ contains $e$ for any $x\in X$. But this is clear:
$e=x^{-1}\cdot x$.

Therefore, by Lemma~\ref{jhgfdjdjfd}, it suffices to check that
$\cK\times T\to \bP_{vc}$ is flat. Let $U\subset \bP_{vc}\times T$
be the pullback of the universal family in $\Hilb(\bP)\times\bP$.
The morphism $U\to \bP_{vc}$ is then flat. But there exists  a
natural map
$$T\times \cK\to U, \quad (t,[Y])\mapsto (t\cdot [Y],t),$$
and this is an isomorphism with inverse $([Y],t)\mapsto
(t,t^{-1}[Y])$. In particular the natural morphism $T\times \cK\to
\bP_{vc}$ is flat and isomorphic to $U\to \bP_{vc}$.
\end{proof}

\begin{Remark}
L.~Lafforgue suggested me that there should be a non-constructive
proof of existence of tropical compactifications using an
appropriate generalization of Raynauld's flattening stratification
theorem over a toric stack $\bP/T$.
\end{Remark}

Now we prove Theorem~\ref{MainTh} as an immediate corollary of
previous results.

\begin{proof}[Proof of Theorem~\ref{MainTh}]
Let $\oX_{vc}\subset\bP_{vc}$ be a visible contour compactification.
It is tropical by Theorem~\ref{KJFKJFKD}. By
Proposition~\ref{:KLJGLJFKF}, $\oX'$ is then tropical for any proper
morphism $\bP'\to\bP_{vc}$. We can assume therefore that $(X,\bP)$
is tropical and $\bP$ is smooth by taking a toric resolution of
singularities of $\bP_{vc}$. Since a boundary of a smooth toric
variety has normal crossings, $\oX$ has ``combinatorial normal
crossings'' by dimension reasons (a structure map is faithfully flat
and therefore equidimensional). Finally, assume that $r=\dim X$ and
$p\in\cap B_i$. The $T$-orbit of $\bP$ containing~$p$ is defined by
a regular sequence $f_1,\ldots,f_r$ (being a closed orbit of the
smooth toric variety). Since the structure map is flat, the pull
back of this sequence is a regular sequence $g_1,\ldots,g_r$ on
$\oX\times T$. It is immediate from definitions that this sequence
cuts out $\{p\}\times T\subset\oX\times T$ set-theoretically. Let
$x_1,\ldots,x_n$ be coordinates on $T$ and let
$t=(\lambda_1,\ldots,\lambda_n)\in T$ be a point such that $(p,t)$
does not belong to the union of embedded components of the scheme
defined by the ideal $(g_1,\ldots,g_r)$. Then the sequence
$(g_1,\ldots,g_r,x_1-\lambda_1,\ldots,x_n-\lambda_n)$ is regular on
$X\times T$. Passing to the localization and changing the order in
the regular sequence we see that $(g_1(t),\ldots,g_n(t))$ is a
regular sequence in the maximal ideal of $p$ on $\oX$.
\end{proof}

\begin{Definition}
Let $\bP$ be an auxiliary compactification of $T$. In the notation
of Theorem~\ref{KJFKJFKD}, the complete fan $\cF_{Gr}$ of $\bP_{Gr}$
is called the {\em Groebner fan} of $(X,\bP)$. The {\em tropical
fan} $\cF_{vc}$ of $\bP_{vc}$ is a subfan of $\cF_{Gr}$ of dimension
$\dim X$. It follows from Proposition~\ref{:KLJGLJFKF} and
Theorem~\ref{KJFKJFKD} that $\cF_{vc}$ is supported on~$\cA$.
\end{Definition}

If an auxiliary toric variety $\bP$ is a projective space,
$\cF_{vc}$ is a {\em tropical variety} of Sturmfels and Speyer
\cite{SS}, though not quite, as we explain below.

\begin{Review}\textsc{Initial Ideals \cite{SS}.}\label{InitialIdeals}
If $\bP$ is a projective space then there is an alternative
description of a tropical fan provided by initial ideals. Let
$\pi:\,V\setminus\{0\}\to\bP$ be the canonical projection. It is
convenient to temporarily change the notation related to tori. The
torus $T$ will be denoted by $T'$ and $T$ will denote
$\pi^{-1}(T')\subset V$. Then $T'=T/k^*$ and we can consider all
introduced $T'$-toric varieties as $T$-varieties and use the
notation $\Lambda$, $\cF_{trop}$, $\cF_{Gr}$, $\ldots$ for
$T$-objects. Let $I\subset k[V]$ be the homogeneous ideal of $\oX$.
Let $\Lambda^\vee_+\subset\Lambda^\vee$ be the standard octant (the
fan of a toric variety $V$).

Any $\chi\in\Lambda^\vee_+$ gives rise to a grading
$$k[V]=\bigoplus_{n\in\bZ}k[V]_n,
\quad k[V]_n=\{f\in k[V]\,|\,\chi(z)\cdot f=z^nf\}.$$ For any $f\in
k[V]$, we define the initial term $in_\chi(f)$ as the homogeneous
component of $f$ of the minimal possible degree. We define the {\em
initial ideal}
$$I_\chi=\{in(f)\ |\ f\in I\}.$$
Then it is well-known that $I_\chi$ describes the flat degeneration,
i.e.~$I_\chi$ is a homogeneous ideal of $\oX_\chi$, where
$$[\oX_\chi]=\lim_{z\to0}\chi(z)[\oX]$$
is the flat limit. In particular, we can describe the Groebner fan
as follows: $\chi',\chi''\in\Lambda^\vee_+$ belong to the interior
of the same cone of $\cF_{Gr}$ if and only if
$I_{\chi'}=I_{\chi''}$. A fan  $\tilde\cF_{Gr}$ defined this way is
actually a refinement of $\cF_{Gr}$ because different initial ideals
can produce the same limit (due to the presence of embedded
components supported on $\{0\}\subset V$). The corresponding
refinement $\cF_{trop}$ of $\cF_{vc}$ is then a subfan of
$\tilde\cF_{Gr}$ described as follows: $\chi\in\cF_{trop}$ if and
only if $I_\chi$ contains no monomials \cite{SS}.
\end{Review}

\section{Very Affine Varieties}

It is easy to see that a connected affine variety $X$ is a (closed)
subvariety of a torus if and only if the ring of regular functions
$\cO(X)$ is generated by units $\cO^*(X)$. In this case we call $X$
{\em very affine} following a suggestion of Kapranov.

It follows from the definition that any very affine variety $X$ is a
subvariety of an {\em intrinsic torus} defined as $T=\Hom(\Lambda,
k^*)$, where $\Lambda=\cO^*(X)/k^*$ (it is well-known that it is a
finitely generated free $\bZ$-module). The embedding $X\subset T$ is
unique up to a translation by an element of $T$. In our main
examples (complements to hyperplane arrangements and open cells of
Grassmannians) that we study in subsequent sections, the tori $T$
are, in fact, intrinsic.

Any morphism $X\to X'$ of very affine varieties induces
the homomorphism $\cO^*(X')\to\cO^*(X)$ by pull-back and therefore
induces the homomorphism of intrinsic tori $T\to T'$.
It also obviously induces the map of amoebas $\cA\to\cA'$.

\begin{Proposition}\label{surj}
If $X\to X'$ is a dominant morphism of very affine varieties
(e.g.~an open immersion) then the corresponding map of amoebas is surjective.
\end{Proposition}

\begin{proof}
We consider $X$ and $X'$ as closed subvarieties of their intrinsic
tori $T$ and~$T'$. Let $f:\,X\to X'$ be a dominant morphism. We
denote the induced map of amoebas by the same letter
$f:\,\cA\to\cA'$. Let $Z$ be a connected Weyl divisor of $X'$ that
contains $X'\setminus f(X)$. We consider $Z$ as a subvariety of the
torus $T'$, let $\cB$ be the corresponding amoeba in
$(\Lambda'_\bQ)^\vee$. The crucial observation \cite{EKL} is that
$\cB$ (or more precisely any fan supported on it) is equidimensional
of dimension equal to the Krull dimension of $Z$ while the dimension
of $\cA'$ is equal to $\dim X'$.

It follows that for any $p\in\cA'$ there exists $x'\in (X'\setminus
Z)(\oK)$ such that $p=\deg x'$. It is clear that there exists $x\in
X(\oK)$ such that $f(x)=x'$. Therefore, $p\in f(\cA)$.
\end{proof}

\begin{Remark}\label{surjectiveamoeba}
This proposition can be used to construct tropical compactifications
of $X'$: find a ``simple'' variety $X$ that dominates $X'$, find an
amoeba $\cA$ of $X$, find an amoeba $\cA'$ of $X'$ as the image of
$\cA$. Then find a fan supported on $\cA'$ and try to prove
(e.g.~using local calculations in toric charts) that the
corresponding compactification is tropical. Applications of this
algorithm to the study of moduli of Del Pezzo surfaces will be
published elsewhere.
\end{Remark}

\begin{Remark}\label{ressing}
It follows (using Proposition~\ref{:KLJGLJFKF}) that for any
tropical compactification $\oX'$, there exists a tropical
compactification $\oX$ and a morphism $\oX\to\oX'$ that extends
$X\to X'$. This motivates the following conjectures: any algebraic
variety $X$ contains a sch\"on very affine open subset $U$
(resp.~any algebraic variety $X$ is dominated by a sch\"on very
affine open subset $U$). This conjecture is interesting because it
obviously implies the resolution of singularities (resp.~the weak
resolution of singularities) and might be equivalent to it.
\end{Remark}

\section{Complements of Hyperplane Arrangements}\label{Hyperplane}
Consider the complement of a hyperplane arrangement
\begin{equation}
X=\bP^r\setminus\{H_1,\ldots,H_n\}.\tag{$\circ$}
\end{equation}
We assume that the arrangement is essential and connected.

$\bP^r$ admits an embedding in $\bP^{n-1}$ such that the hyperplanes
$H_i$ are intersections with coordinate hyperplanes of~$\bP^{n-1}$.
It follows that $X$ is a closed subvariety of a torus
$T=(k^*)^{n-1}$. It is easy to see that this torus is intrinsic.

We take $\bP^{n-1}$ as the auxiliary compactification of the visible
contour construction of Definition~\ref{defvisiblecontour}.
Connectedness of the arrangement implies $T_X=\{e\}$. The component
of the Hilbert scheme of $\bP^{n-1}$ containing $\oX=\bP^r$ is
obviously the Grassmannian $G(r+1,n)$. Let $[\oX]\in G(r+1,n)$ be
the corresponding point. The toric Groebner variety $\bP_{Gr}$ is
the normalization of a closure of $T\cdot[\oX]\subset G(r+1,n)$.

It follows that $\oX_{vc}$ is equal to the closure of $X^{-1}\cdot
[\oX]$ in $\bP_{Gr}$ (in this strange looking formula $X^{-1}$ means
inverses of elements of $X\subset T$). This is precisely Kapranov's
visible contour (he takes the closure in $G(r+1,n)$ rather than
$\bP_{Gr}$ but this does not matter: $\oX_{vc}\subset\bP_{vc}$ and
$\bP_{vc}$ is isomorphic to its image in $G(r+1,n)$ (i.e.~this image
is normal) by a result of White \cite{W}).

We continue with the proof of Theorem~\ref{:LKJHGKGKJYF}.

Let $\Psi:\,T\times\oX_{vc}\to\bP_{vc}$ be the multiplication map.
Since the universal family of the Grassmannian is smooth (being a
projective bundle), we see as in the proof of Theorem~\ref{KJFKJFKD}
that $\Psi$ is not just flat but smooth, i.e.~$X$ is sch\"on. By
Theorem~\ref{KLHGKFJGFhgghhg}, $\oX_{vc}$ is normal and has toroidal
singularities. Let $G_e\subset G(r+1,n)$ be the subvariety of
subspaces containing $e=(1,\ldots,1)\in T$. Obviously, $G_e$ is
isomorphic to $G(r,n-1)$ and is regularly embedded in $G(r+1,n)$
with normal bundle being the universal quotient bundle of~$G_e$. By
Theorem~\ref{KJFKJFKD}, $\oX_{vc}$ is a scheme-theoretic
intersection of $\bP_{vc}$ with $G_e$. In particular, the normal
bundle of $\oX_{vc}\subset\bP_{vc}$ in this case is the pullback of
the universal quotient bundle of~$G_e$. By
Theorem~\ref{KLHGKFJGFhgghhg}, the log canonical bundle of
$\oX_{vc}$ is the determinant of its normal bundle. It follows that
it is very ample being the restriction of a Pl\"ucker polarization
of $G_e$. This finished the proof of Theorem~\ref{:LKJHGKGKJYF}.
\qed\medskip

\begin{Example}\textsc{\cite{CQ2}.}\label{LKGHKFKJGFGFhh}
Assume that $X=\bP^2\setminus\{L_1,\ldots,L_n\}$. Let the {\em
multiplicity} of $p\in\bP^2$ be the number of lines passing through
it. Let $\oX_{wond}$ be the blow up of $\bP^2$ in the set of points
of multiplicity at least $3$. Restricting the logcanonical line
bundle on (proper transforms of) lines and exceptional divisors, we
see that it is ample, and therefore $\oX_{wond}=\oX_{vc}$ with the
following exception. Suppose the configuration of lines contains a
line~$L$ and two points $a,b\in L$ such that any other line passes
through $a$ or~$b$. In this case the restriction of the logcanonical
line bundle on $L$ is trivial, and therefore
$\oX_{vc}\simeq\bP^1\times\bP^1$ (the blow-down of the proper
transform of $L$ in $\oX_{wond}=\Bl_{a,b}\bP^2$).
\end{Example}

This example generalizes to the higher dimensions as follows. Since
$\oX_{vc}$ is the log canonical model of $X$, it is natural to
wonder how $\oX_{vc}$ is related to compactifications of $X$ with a
normally crossing boundary. They can be constructed as subsequent
blow ups of $\bP^r$ along (proper transforms of) projective
subspaces required to make the proper transform of the hyperplane
arrangement a divisor with normal crossings. For example, one can
blow up all possible partial intersections of hyperplanes and obtain
the variety $BL$ of \cite[\S5]{CQ2}. The most economical blow up is
the so-called {\em wonderful compactification} $\oX_{wond}$
constructed in \cite{CP}. In fact, there is a natural poset of
compactifications with normal crossings constructed in~\cite{FY}
such that $BL$ is the maximal and $\oX_{wond}$ is the minimal
element of this poset. These compactifications $\oX_{nest}$ are
indexed by the so called nesting sets. Moreover, any $\oX_{nest}$ is
the closure of $X$ in a certain smooth toric variety $\bP_{nest}$.
It is shown in \cite{FS} by a combinatorial argument that the fan
$\cF_{nest}$ refines $\cF_{vc}$ and in particular $\oX_{nest}$ is
tropical by Proposition~\ref{:KLJGLJFKF} and there is a canonical
map $\oX_{nest}\to\oX_{vc}$. This map is log crepant by
Theorem~\ref{KLHGKFJGFhgghhg}, i.e.~the logcanonical divisor $K+B$
is globally generated on $\oX_{nest}$ and gives a regular map
$\oX_{nest}\to\oX_{vc}$ identical on $X$.
Example~\ref{LKGHKFKJGFGFhh} shows that $\oX_{wond}$ is not
necessarily equal to $\oX_{vc}$. The combinatorial condition to
guarantee $\oX_{wond}=\oX_{vc}$ was found in \cite{FS}
(geometrically it means that any strata of $\oX_{wond}$ is of log
general type).

\section{Tropical Quotients of Grassmannians}\label{tropicalquotients}

\begin{Review}\textsc{Tropical Compactifications as Compact Quotient Spaces.}\label{JJTHGFHGHHHGHG}
Let $X$ be a closed subvariety of a torus $T$. Here we will be
interested in the case $T_X\ne\{e\}$. Since $\oX_{trop}$
compactifies $X/T_X$, it is natural to compare $\oX_{trop}$ with
other known compact quotients. The most closely related ones are
Hilbert/Chow quotients, see \cite{KSZ,CQ1}. The construction goes as
follows. We first take any auxilliary compact toric variety $\bP$ of
$T$ as in the construction of $\oX_{vc}$. Let $Z\subset \oX$ be the
closure of a generic $T_X$-orbit in $\oX$. Then we have an embedding
$$T/T_X\hookrightarrow\Hilb(\bP),\quad t\mapsto t\cdot [Z],$$
where we take a component of a Hilbert scheme containing $[Z]$ (note
that the construction of $\oX_{vc}$ uses a completely different
component). Let $\bP_{Hilb}$ be the normalization of the closure of
$T/T_X$ in $\Hilb(\bP)$. The closure of $X/T_X$ in $\bP_{Hilb}$ is
the {\em Hilbert quotient}  $\oX_{Hilb}$ (one can also use a Chow
variety of $\bP$ or its multigraded Hilbert scheme to obtain closely
related quotients).

The difference between $\oX_{Hilb}$ and $\oX_{vc}$ can be summarized
as follows: $\oX_{Hilb}$ parametrizes $T$-translates of $Z$ (and
their limits -- broken toric varieties) in $\bP$ that are contained
in $\oX$\footnote{More precisely, let $S$ be the set of these
generic orbits and their limits. $S$ has a natural scheme structure.
In general $S$ is reducible (example is given in \cite{CQ2}) and
$\oX_{Hilb}$ is the reduction of the principal component of $S$},
while $\oX_{trop}$ parametrizes $T$-translates of $\oX$ (and their
limit degenerations) that contain $Z$. The reason why $\oX_{vc}$ is
better than $\oX_{Hilb}$ is because inclusion-schemes (subfamily of
a family of subschemes contained in a fixed subvariety) in general
behave worse than containment-schemes (subfamily of a family of
subschemes that contain a fixed subvariety)

\begin{Definition}
Let $\bP_{Laf}\subset\bP_{Hilb}$ be the open chart (introduced by
Lafforgue in the case of the Grassmannian) consisting of $T$-orbits
that intersect $\oX_{Hilb}$. Let $\cF_{Laf}$ be the corresponding
fan.
\end{Definition}

$\oX_{Hilb}$ is often not tropical. Here we give a stupid example
with $T_X=e$, an example with $T_X\ne\{e\}$ is given in
Remark~\ref{G(3,9)}. Let $X$ be a complement to a not normally
crossing connected hyperplane arrangement \eqref{KContour}.
Obviously
$$\oX_{Hilb}=\bP^r\subset\bP_{Hilb}=\bP^{n-1}$$
because $T_X=\{e\}$. Since the boundary does not have combinatorial
normal crossings, the structure map is not flat by
Theorem~\ref{MainTh}. On the other hand, $\oX_{vc}$ magically
incorporates all necessary blow ups. The special configuration of
the Example~\ref{LKGHKFKJGFGFhh} also shows that in general there is
no regular map $\oX_{vc}\to\oX_{Hilb}$.
\end{Review}

\medskip

In the remainder of this section we study the quotient of a
Grassmannian by a maximal torus. Let $X=G^0(r,n)$ be the open cell
in the Grassmannian $G(r,n)$ given by non-vanishing of all Pl\"ucker
coordinates. Consider the Pl\"ucker embedding
$$X=G(r,n)\subset\bP=\bP(\Lambda^rk^n).$$
The complement to coordinate hyperplanes in $\bP$ is a torus $T$
that contains $X$ as a closed subvariety. It is easy to see that $T$
is an intrinsic torus of $X$. $T_X=(k^*)^n/k^*$ is the diagonal
torus. We fix $\bP$ as an auxiliary toric variety necessary to
construct the Hilbert quotient $\oX_{Hilb}$ and the generalized
visible contour compactification $\oX_{vc}$ and compare them. In
fact, instead of $\oX_{vc}$, we will consider its slight
modification $\oX_{trop}$ defined in \ref{InitialIdeals} using
initial ideals. I don't know if $\oX_{vc}=\oX_{trop}$. The fan
$\cF_{trop}$ is the {\em tropical Grassmannian} of Sturmfels and
Speyer \cite{SS}. $\oX_{Hilb}$ is Kapranov's Chow quotient of
Grassmannian introduced in \cite{CQ1}.

\begin{Remark}\label{G(3,9)}
By the Gelfand--Macpherson correspondence $X/T_X$ is identified with
the moduli space $X(r,n)$ of arrangements of $n$ hyperplanes in
$\bP^{r-1}$ in linearly general position. By \cite{HKT},
$\oX_{Hilb}$ has a functorial description as the moduli space of
stable pairs of log general type (and thus it serves as a template
for the study of moduli spaces of varieties of general type).
Description of fibers of the universal family over $\oX_{Hilb}$
(higher dimensional analogs of stable rational curves) were found in
\cite{Laf}, their natural crepant resolutions were found in
\cite{CQ2}. However, by \cite{CQ2} the geometry of $\oX_{Hilb}$ is
terrible: for any scheme $S$ of finite type over $\Spec\bZ$, there
is a strata in $\oX_{Hilb}$ isomorphic to the open subset $U\subset
S\times\bA^r$ such that the projection $U\to S$ is onto. Also, if
$n\ge9$ then $\oX_{Hilb}$ is not a log canonical model of $X$.
Moreover, it is proved in \cite{CQ2} that the fan of $\bP_{Hilb}$
for $n\ge9$ contains cones of dimension bigger than dimension of the
amoeba. In particular, $\oX_{Hilb}$ is not tropical for $n\ge9$.
\end{Remark}

\begin{Theorem}\label{JHGjjjFGJHG}
There exists a toric morphism $\bP_{trop}\to\bP_{Laf}$.
\end{Theorem}

In particular, there exists a morphism $\oX_{trop}\to\oX_{Hilb}$. In
some cases we can say more. The following theorem is a simple
compilation of known results:

\begin{Theorem}\label{m0n}
Let $X=G^0(2,n)$. Then $\oX_{Hilb}=\oX_{trop}=\oM_{0,n}$, the
Grothendieck--Knudsen moduli space of stable rational curves.
\end{Theorem}

\begin{proof}
Fans $\cF_{Laf}$ and $\cF_{vc}$ were calculated in \cite{CQ1} and
\cite{SS}, respectively. They are the same, the so-called {\em space
of phylogenetic trees}. Therefore
$$\oX_{trop}=\oX_{Hilb}=\oM_{0,n}$$
the last equality is a result of Kapranov \cite{CQ1}.
\end{proof}

\begin{Remark}
Notice that $M_{0,n}$ is also the complement of the {\em braid
arrangement}
$$\{x_i=x_j\,|\,1\le i<j\le n-1\}\subset \bP^{n-3},$$
where $\sum x_i=0$. In particular, it has another visible contour
compactification defined in the framework of
Theorem~\ref{LKJGKFKJGFBVCH}. Since $\oM_{0,n}$ is the log canonical
model of $M_{0,n}$, this latter compactification is also isomorphic
to $\oM_{0,n}$ by Theorem~\ref{LKJGKFKJGFBVCH} and gives its
Kapranov's blow up model defined in \cite{CQ1}.
\end{Remark}

The next result is joint with Hacking and Keel. Its proof will
appear elsewhere. Conjecturally, its analogue holds for $G(3,7)$ and
$G(3,8)$.

\begin{Theorem}\label{KHGJjhj}
Let $X=G^0(3,6)$. Then $X$ is sch\"on, $\oX_{Hilb}$ is its log
canonical model. It has $40$ isolated singularities. $\oX_{trop}$ is
a crepant resolution of $\oX_{Hilb}$ .
\end{Theorem}

The rest of the paper is occupied by the proof of
Theorem~\ref{JHGjjjFGJHG}.

\begin{proof}
Let $\Delta=\Delta(r,n)\subset\bR^n$ be the {\em hypersymplex},
i.e.~the convex hull of points
$$e_I:=\sum_{i\in I}e_i,\quad\hbox{\rm where}\quad I\subset N:=\{1,\ldots,n\},\quad |I|=r.$$
It is well-known that all these points are vertices of $\Delta$ and we frequently identify
$\Delta$ with its set of vertices if it does not cause confusion.

\begin{Definition}
Let
$$\cH=\{H_1,\ldots,H_n\}\subset\bP^{r-1}$$
be an essential hyperplane arrangement (not all hyperplanes pass
through a point). In contrast to the Section~\ref{Hyperplane}, we
allow some hyperplanes to appear more than once (to be multiple).
Then the {\em matroid polytope} $P_\cH\subset\Delta(r,n)$ of $\cH$
is the convex hull of vertices $e_I$ such that
$H_{i_1},\ldots,H_{i_r}$ are linearly independent.
\end{Definition}

\begin{Remark}
So for example $\Delta$ itself is a matroid polytope of a configuration
of hyperplanes in linearly general position.
$P_\cH$ has the maximal dimension if and only if  $\cH$
is a connected configuration, i.e.~has no projective automorphisms.
\end{Remark}

\begin{Definition}
For any $i_0\in N$, we define a {\em restricted configuration} $R_{i_0}(\cH)\subset H_{i_0}$
as the set $H_i\cap H_{i_0}$ for any $i$ such that $H_i\ne H_{i_0}$
(so the indexing set of $R_{i_0}(\cH)$ can be strictly smaller than $N\setminus\{i_0\}$).
For any subset $K\subset N$,
we define a {\em contracted configuration} $C_K(\cH):=\cH\setminus\{H_k\}_{k\in K}\subset\bP^{r-1}$.
We say that $\cH$ is {\em more constrained} than $\cH'$ if any linearly independent subset
of~$\cH$ is linearly independent in $\cH'$ but not vice versa (this is equivalent
to the strict inclusion $P_{\cH}\subset P_{\cH'}$).
More generally, if $\cH$ is indexed by~$M$, $\cH'$ is indexed by~$N$, and $M\subset N$
then we say that $\cH$ is more constrained than $\cH'$ if $\cH$
is more constrained than $C_{N\setminus M}(\cH')$.
\end{Definition}

\begin{Lemma}\label{LKJGKKGHGJFKJjj}
Let $r>2$. Let $\cH,\cH'\subset\bP^{r-1}$ be connected configurations of $n$ hyperplanes
such that $\cH$ is more constrained than $\cH'$. Then there exists $i_0\in N$
such that $R_{i_0}(\cH)$ is connected and more constrained than $R_{i_0}(\cH')$
(which in this case is of course also automatically connected).

Let $\cH,\cH'\subset\bP^1$ be connected configurations of $n>4$
points such that $\cH$ is more constrained than $\cH'$. Then there
exists an index $i_0\subset N$ such that $C_{i_0}(\cH)$  is
connected and more constrained than $C_{i_0}(\cH')$ (which in this
case is also connected).
\end{Lemma}

\begin{proof}
Let $\{H_1,\ldots,H_k\}\subset\cH$ be the minimal set of linearly dependent hyperplanes
such that the corresponding hyperplanes in $\cH'$ are linearly independent.
Then we can find linear equations of these hyperplanes such that
$f_1,\ldots,f_{k-1}$ are linearly independent and $f_k=f_1+\ldots+f_{k-1}$.
Let $\hat H_k\in\cH$ be a hyperplane such that $\{H_1,\ldots,H_{k-1},\hat H_k\}$
is linearly independent. If $k=2$ then we assume, in addition, that
$H_1',H_2',\hat H_2'$ is linearly independent (here we use $r>2$).
Finally, we construct a basis
$\{H_1,\ldots,H_{k-1},\hat H_k,H_{k+1},\ldots,H_r\}\subset\cH$.
We define a graph $\Gamma$ with vertices $\{1,\ldots,r\}$ as follows:
$i$ and $j$ are connected by an edge if and only if there exists a hyperplane $H\in\cH$
such that in the expression of its equation in our basis $i$-th and $j$-th
coefficients are both nontrivial. Then it is immediate from the calculation
of a (trivial) automorphism group of a configuration in our basis that this graph is connected.
Therefore, it has at least two vertices that can be dropped without
ruining the connectedness. We take $i_0$ to be any of these vertices unless $k=2$
in which case we take $i_0$ to be any of these vertices not equal to $1$.
Now all claims follow by immediate inspection.

The statement about $\bP^1$ is obvious.
\end{proof}

\begin{Review}\textsc{Pictures.}
Let $R\subset K$ be the DVR of power series.
Let
$$\cH(z)=\{H_1(z),\ldots,H_n(z)\}\subset\bP^{r-1}(K)$$
be the collection of (one-parameter families) of hyperplanes in
linearly general position. Any choice of a $K$-frame $F$ (the
projectivization of a basis in $K^r$) gives a canonical way to
extend this $K$-point of $\bP^{r-1}$ to the $R$-point. Notice that
frames related by an element of $PGL_r(R)$ will give isomorphic
families over $\Spec R$, i.e.~points of the homogeneous space
$PGL_r(K)/PGL_r(R)$ (the so called {\em affine building}) give a
family defined upto an isomorphism. The special fiber of this family
is $\bP^{r-1}$ with a limiting configuration of hyperplanes
$\cH^F(0)=\{H^F_1(0),\ldots,H^F_n(0)\}\subset\bP^{r-1}$ that depends
on the choice of a frame. Most of the time this configuration is not
even essential, frames that give essential configurations is, by
definition, a {\em membrane} of \cite{CQ2}.\footnote{A membrane is a
polyhedral complex surprisingly homeomorphic to the nonarchimedean
amoeba of a very affine variety $\bP^{r-1}(K)\setminus\cH(z)$
defined over $K$ \cite{CQ2}.} It is proved there that matroid
polytopes of these essential configurations give a {\em matroid
decomposition}, i.e.~a paving of $\Delta$ by a finite number of
matroid polytopes. As mentioned above, cells (polytopes of the
maximal dimension) of this decomposition correspond to limiting
configurations that are connected.
\end{Review}

\begin{Definition}
We say that a configuration of hyperplanes $\cH(z)$
over $K$ in linearly general position is {\em more constrained} than $\cH'(z)$
if the matroid decomposition induced by $\cH'(z)$ is coarser
than the matroid decomposition induced by $\cH(z)$.
\end{Definition}

\begin{Lemma}\label{lkjghlGJGJGKJKJ}
Let $r>2$. Let $\cH(z),\cH'(z)\subset\bP^{r-1}(K)$ be configurations
of $n$ hyperplanes in linearly general position such that $\cH(z)$
is more constrained than $\cH'(z)$. Then there exists $i_0\in N$
such that $R_{i_0}(\cH(z))$ is more constrained than
$R_{i_0}(\cH'(z))$ (also notice that $R_{i_0}(\cH(z))$ and
$R_{i_0}(\cH'(z))$ are automatically in linearly general position
because $\cH(z)$ and $\cH'(z)$ are).

Let $\cH(z),\cH'(z)\subset\bP^1(K)$ be configurations of different
$n>4$ points such that $\cH(z)$ is more constrained than $\cH'(z)$.
Then there exists an index $i_0\subset N$ such that $C_{i_0}(\cH)$
is more constrained than $C_{i_0}(\cH')$ (also notice that
$C_{i_0}(\cH(z))$ and $C_{i_0}(\cH'(z))$ are automatically in
linearly general position because $\cH(z)$ and $\cH'(z)$ are).
\end{Lemma}

\begin{proof}
There are frames $F$ and $F'$ for $\cH(z)$ and $\cH'(z)$ that
produce matroid polytopes of maximal dimension embedded in each
other. This means that the corresponding limiting configuration
$\cH^F(0)$ is more constrained then $\cH^{F'}(0)$. Now apply
Lemma~\ref{LKJGKKGHGJFKJjj}.
\end{proof}

Let $X=G^0(r,n)$. We identify the Pl\"ucker vector space $V$
with~$k^\Delta$, the intrinsic torus $T$ with $(k^*)^\Delta$ (using
the convention about tori from \eqref{InitialIdeals}), and the
lattice $\Lambda^\vee$ with $\bZ^\Delta$. Let $\cI$ be the ideal of
the cone over $G(r,n)$ in $k^\Delta$.

\begin{Review}\textsc{Description of Fans.}\label{fandescr}
As a set, the ``tropical Grassmannian'' $\cF_{trop}$ is described in
\cite[Theorem 3.8]{SS}. Namely, $\omega\in\bZ^\Delta$ belongs to
$\cF_{trop}$ if and only if there exists $[L]\in X_{K}$ (or
equivalently the one-parameter family of configurations of
hyperplanes $\cH(z)\subset\bP^{r-1}(K)$ as above by the
Gelfand--Macpherson transform) such that
\begin{equation}\label{KJHLKHGLKHGKHG}
\omega(I)=\deg\Delta_I(L),\cooltag
\end{equation}
where $\Delta_I$ is the Pl\"ucker coordinate that corresponds to
$I\subset N$, $|I|=r$, and $\deg:\,K\to\bZ$ is the standard
valuation. The matroid decomposition corresponding to $\cH(z)$ is
can be recovered from $\omega$ as follows. Let
$\Delta_\omega\subset\Delta\times\bZ\subset\bR^{n+1}$ be the convex
hull of points $(I,\omega(I))$ for $I\subset N$, $|I|=k$. Then all
these points are vertices of $\Delta_\omega$ and projections of
faces of $\Delta_\omega$ on $\Delta$ give a matroid decomposition of
$\Delta$. By \cite{Laf} or \cite{CQ1}, two functions in $\bZ^\Delta$
belong to the interior of the same cone in $\cF_{Laf}$ if and only
if they induce the same matroid decomposition and the inclusion of
cones in $\cF_{Laf}$ corresponds to the coarsening of the paving.
\end{Review}

\begin{Review}\textsc{Face Maps.}
We have to remind definitions of face maps and cross-ratios
\cite{CQ1,Laf,CQ2}. For any $i_0\in N$, restriction and contraction
of configurations of hyperplanes in linearly general position can be
extended to morphisms
$$\oX_{Hilb}(r,n)\to \oX_{Hilb}(r-1,n-1),\quad
\oX_{Hilb}(r,n)\to \oX_{Hilb}(r-1,n).$$
These morphisms are induced by morphisms of toric varieties
related to
the natural inclusions
$$\Delta(r-1,n-1)=\{p\in\Delta(r,n)\ |\ x_{i_0}(p)=1\}$$
and
$$\Delta(r-1,n)=\{p\in\Delta(r,n)\ |\ x_{i_0}(p)=0\}.$$
More precisely, let us interpret $\cF_{Laf}\subset\bQ^\Delta$ as a set
of functions on the hypersymplex.
Then restriction and contraction are toric morphisms
given by maps of fans obtained by
restricting functions to the boundary of the hypersymplex.
\end{Review}

In particular, for any collection of hyperplanes
$$\{H_{i_1},\ldots,H_{i_4};H_{j_1},\ldots,H_{j_{r-2}}\},\quad\hbox{\rm where}\quad
I\cap J=\emptyset,$$
there is a {\em cross-ratio map}
$$\oX_{Hilb}\to \oX_{Hilb}(2,4)=\oM_{0,4}=\bP^1:$$
contract all hyperplanes not in $I\cup J$ and then consecutively
restrict on all hyperplanes from $J$. The cross-ratio map is induced
by the toric morphism $\bP_{Laf}\to\bP^2\setminus\{e_1,e_2,e_3\}$
given by the inclusion of the octahedron
$\Delta(2,4)\subset\Delta(r,n)$, where $\Delta(2,4)$ has vertices
$$e_{i_a}+e_{i_b}+e_{j_1}+\ldots+e_{j_{r-2}}\quad\hbox{\rm for}\quad
(a,b)\subset\{1,2,3,4\}.$$

One can predict whether the cross-ratio  is equal to $0$, $1$, or
$\infty$ by looking at the decomposition of the octahedron
$\Delta(2,4)$ induced by the matroid decomposition of $\Delta$
because $\{0,1,\infty\}\subset\bP^1$ is the intersection of $\bP^1$
with the divisorial boundary of $\bP^2\setminus\{e_1,e_2,e_3\}$. The
values of $0$, $1$, and $\infty$ correspond to three different ways
to decompose the octahedron as a union of two pyramids.

\begin{Lemma}\label{KJHKJHKHGHG}
Let $\omega\in\cF_{trop}$. Then the initial ideal $\cI_\omega$
determines uniquely whether any cross-ratio is equal to $0$, $1$,
$\infty$, or none of the above.
\end{Lemma}

\begin{proof}
Indeed, $\cI$ includes the simplest Pl\"ucker relation
$$f=\Delta_{12J}\Delta_{34J}-\Delta_{13J}\Delta_{24J}+\Delta_{14J}\Delta_{23J}.$$
This is the unique (up to a scalar) element of $\cI$ of its
$(k^*)^n$ weight (this easily follows by inspection by looking at
other Pl\"ucker relations). Since $\cI_\omega$ does not contain
monomials, there are only $4$ possibilities for $in_\omega(f)$:
either $in_\omega(f)=f$, i.e.
$$\omega_{12J}+\omega_{34J}=\omega_{13J}+\omega_{24J}=\omega_{14J}+\omega_{23J},$$
in which case $\Delta(2,4)$ is obviously not split.
Or $in_\omega(f)$ can be a sum of two monomials, e.g.
$$\Delta_{12J}\Delta_{34J}-\Delta_{13J}\Delta_{24J},$$
in which case
$$\omega_{12J}+\omega_{34J}=\omega_{13J}+\omega_{24J}<\omega_{14J}+\omega_{23J}$$
and $\Delta(2,4)$ is split into $2$ pyramids. In any case we see
that $\cI_\omega$ determines uniquely whether any cross-ratio is
equal to $0$, $1$, $\infty$, or none of the above.
\end{proof}

Now we can prove the Theorem.  We have to show that the interior of
any cone in $\cF_{trop}$ belongs to the interior of a cone in
$\cF_{Laf}$. It suffices to show that if
$\omega,\omega'\in\cF_{trop}\cap\bZ^\Delta$ and $\omega'$ gives a
coarser matroid decomposition than $\omega$ then $\cI_\omega\ne
\cI_{\omega'}$. Choose $\cH(z)$ and $\cH'(z)$ that correspond to
$\omega$ and $\omega'$ as in \ref{fandescr}. Then $\cH(z)$ and
$\cH'(z)$ is more constrained than $\cH'(z)$ by definition.

Applying Lemma~\ref{lkjghlGJGJGKJKJ} consecutively, we can find
subsets $I$ and $J$ with $|I|=4$ and $|J|=r-2$ such that
$R_JC_{(I\cup J)^c}(\cH(z))$ is more constrained than $R_JC_{(I\cup
J)^c}(\cH'(z))$. This means that the octahedron $\Delta(2,4)$
corresponding to $I$ and $J$ is split by $\omega$ and not split by
$\omega'$. By Lemma~\ref{KJHKJHKHGHG}, it follows that
$\cI_\omega\ne \cI_{\omega'}$.
\end{proof}

\end{document}